\documentclass[11pt]{article}
\usepackage{amsmath, amssymb}
\usepackage{array,tabularx,longtable}
\usepackage{graphicx}

\def \qed{\hfill$\Box$}
\renewcommand{\abstractname}{Abstract.}
\renewcommand\abstract{\hfil\break\topsep=0pt\partopsep=0pt\parsep=0pt\itemsep=0pt\relax
\trivlist\item[\hskip\labelsep
{\bfseries\abstractname}]\if!\abstractname!\hskip-\labelsep\fi}

\newcommand{\email}[1]{{(e-mail: #1)}}

\newcommand{\e}{\mathbb{E}}
\newcommand{\p}{\mathbb{P}}
\newcommand{\f}{\mathcal{F}}

\newcommand{\xn}{X^N}
\newcommand{\xnm}{X^{N,m}}
\newcommand{\xnk}{X^N_{t_k}}
\newcommand{\xntm}{\tilde{X}^{N,m}}

\newcommand{\dw}{\Delta W}
\newcommand{\dwt}{\Delta \tilde{W}}

\newcommand{\yn}{Y^N}
\newcommand{\ynr}{Y^{N,R}}

\newcommand{\pynr}{y^{N,R}}
\newcommand{\pynrm}{y^{N,R,M}}

\newcommand{\zn}{Z^N}
\newcommand{\znr}{Z^{N,R}}

\newcommand{\pznr}{z^{N,R}}
\newcommand{\pznrm}{z^{N,R,M}}

\newcommand{\mc}[1]{\mathcal{#1}}
\renewcommand{\r}{\mathbb{R}}
\newcommand{\etal}{{\em et al. }}
\newcommand \A[1]{{\bf (#1)}}
\newcommand{\1}{\mathbf{1}}
\newcommand{\id}{\mathrm{Id}}
\newcommand{\nn}{\nonumber}

\newtheorem{theorem}{Theorem}

\def \dd{{\rm d}}

\title{Numerical simulation of BSDEs using empirical regression methods: theory and practice\footnote{This article has been presented at the Fifth Colloquium on BSDEs at Shangai on 29th May - 1st June 2005. E. Gobet thanks the organizers of this wonderful meeting, especially professors S. Peng and S. Tang for having invited him.  This work has been financially supported by  {\em Association Nationale de la Recherche Technique}, {\em \'Ecole Polytechnique} and {\em \'Electricit\'e de France}.}}

\date{November 4, 2005}
\author{Emmanuel Gobet\footnote{INP Grenoble - IMAG-LMC - BP 53 - 38041 Grenoble Cedex - FRANCE
\email{emmanuel.gobet@imag.fr}. Corresponding author.}, Jean-Philippe Lemor\footnote{\email{jp.lemor@libertysurf.fr}}}

\begin{document}
\maketitle

\parindent=0 em

\begin{abstract} This article deals with the numerical resolution of backward stochastic differential equations. Firstly, we consider a rather general case where the filtration is generated by a Brownian motion and a Poisson random measure. We provide a simulation algorithm based on iterative regressions on function bases, which coefficients are evaluated using Monte Carlo simulations. We state fully explicit error bounds. Secondly, restricting to the case of a Brownian filtration, we consider reflected BSDEs and adapt the previous algorithm to that situation. The complexity of the algorithm is very competitive and allows us to treat numerical results in dimension 10.
\end{abstract}

\section*{Introduction}
Since the early nineties, there has been an increasing interest for BSDEs. Due to numerous contributions, we now have a better understanding of the features of these equations. This is still a very active field of research as this special volume may show. BSDEs have a wide range of applications, especially in finance, see for instance \cite{elka:peng:quen:97}. Our concern is rather related to the simulation of BSDEs, and regarding these aspects, the contributions in the literature are more recent and less definitive. Actually, the design of efficient algorithms that are able to solve general BSDEs in any reasonable dimension is far to be solved. 

For a review of different contributions, we refer to \cite{gobe:lemo:wari:05}. Here we present a new algorithm. It is designed for rather general BSDEs which are coupled with a forward Markov process $X$. 
 The scheme that we propose has two ingredients. The first one is somewhat standard and consists in approximating the BSDE by a discrete-time backward dynamic programming equation. Then one needs to compute a sequence of regression functions: our choice is to project these functions on a finite-dimensional space spanned by basis functions, which coefficients are computed using a set of simulations of $X$. The details of the algorithm are given in Section \ref{section:1}, in a generalized case where the filtration is generated by a Brownian motion and a Poisson random measure. In that case, $X$ is a jump-diffusion process. We state explicit error estimates, with respect to the parameters of the method which are the time step $h$ (used to derive the discrete-time dynamic programming equation; this parameter may be also used as a time step to simulate $X$), the number of basis functions, the number $M$ of simulations of $X$, the threshold $R$ needed to ensure a stability property (its impact can be usually neglected).

We would like to emphasize several valuable features of our approach. Firstly, the explicit error bound allows us to choose optimally how the above parameters should vary together to achieve a given accuracy. Secondly, unlike the quantization approach \cite{ball:page:03}, here the driver may depend on $Z$. Thirdly, we use only one set of $M$ simulations to evaluate all the regression operators at once, while in \cite{bouc:touz:04} $M\times \#\{\text{discretization times}\}=O(M/h)$ simulations are required. Finally, our scheme makes a very little use of the model for $X$. Indeed, on the one hand, for the implementation one just needs to know how to simulate approximative paths of $X$. On the other hand, for the analysis of regression errors, we use distribution-free techniques (see Gy\"orfi \etal \cite{gyor:kohl:krzy:walk:02}): hence a significant part of error estimates is {\em model-free}. We guess that this is a significant advantage, which allows us to consider general models (for instance, with the Malliavin calculus approach in \cite{bouc:touz:04}, $X$ is assumed to be essentially an elliptic diffusion). 

In Section \ref{section:2}, we consider reflected BSDEs, taking for $X$ a diffusion process. We adapt the previous algorithm into three different variants, depending on how we handle the reflection. It alternatively gives lower and upper solutions.
 
The scheme presented in this article has been first introduced and analyzed in \cite{lemo:gobe:wari:05}: therein, one deals with the framework given in Section \ref{section:1}. We refer the reader to the mentioned reference for the detailed proofs. Materials of Section \ref{section:2} come from the PhD thesis of the second author \cite{lemo:05}.

\section{Generalized backward stochastic differential equation}
 \label{section:1}
We follow the presentation of Barles \etal \cite{barl:buck:pard:97}. Let
$(\Omega,{\cal F},({\cal F}_t)_t,\p)$ be a stochastic basis, where
the filtration satisfies the usual conditions of right-continuity
and completeness. We suppose that the filtration is generated by
the two mutually independent processes: a $\r^q$-valued Brownian
motion $W$ and a Poisson random measure $\mu$ on $\r_+\times E$,
where $E=\r^l\backslash \{0\}$ is equipped with its Borel field
${\cal E}$, with compensator $\nu(\dd t,\dd e)=\dd t \lambda(\dd
e)$, such that $\{\tilde \mu([0,t]\times A)=(\mu-\nu)([0,t]\times
A)\}_{t\geq 0}$ is a martingale for all $A\in {\cal E}$ with
$\lambda(A)<+\infty$. $\lambda$ is assumed to be a $\sigma$-finite
measure on $(E,{\cal E})$ satisfying $\int_E (1\land
|e|^2)\lambda(\dd e)<+\infty$. We consider the $\r^d$-valued
jump-diffusion
\begin{align}
X_t =  x+\int_{0}^{t}b(s,X_s) \dd s + \int_{0}^{t}\sigma(s,X_s)\dd
W_s+\int_{0}^{t}\int_E\beta(s,X_{s^-},e) \tilde \mu(\dd s,\dd e),
\end{align}
which is uniquely defined under the following assumption.
\begin{itemize}
\item [\A{H1}] The functions $b(t,x)$ and $\sigma(t,x)$ are uniformly Lipschitz continuous with respect to $(t,x)\in[0,T]\times \r^d$. For some constant $c$, the function $\beta$ satisfies
$|\beta(t,x,e)|\leq c(1\land |e|)$ and
$|\beta(t,x,e)-\beta(t',x',e)|\leq c(|x-x'|+|t-t'|)(1\land |e|)$
for any $(t,x),(t,x')\in[0,T]\times \r^d$ and $e\in E$.
\end{itemize}
Now consider the following generalized BSDE
\begin{align}
-\dd Y_t = & f(t,X_t,Y_t,Z_t) \dd t - Z_t \dd W_t - \dd L_t,\qquad
Y_T=\phi(X_T) \label{edsr:general},
\end{align}
where $L$ is c\`adl\`ag martingale orthogonal to $W$ (with $L_0=0$). Provided that the driver $f$ is Lipschitz continuous, there is a unique solution $(Y,Z,L)$ (in appropriate $L_2$-spaces, see \cite{barl:buck:pard:97} and \cite{elka:huan:97} for details). To derive explicit error estimates, we impose to the terminal function $\phi$ to be Lipschitz continuous as well. These two regularity assumptions are stated as follows:
\begin{itemize}\item [\A{H2}] For any $(t_1,x_1,y_1,z_1), (t_2,x_2,y_2,z_2)\in [0,T]\times \r^d\times \r \times \r^q$, one has $|f(t_2,x_2,y_2,z_2)-f(t_1,x_1,y_1,z_1)|\leq C_f(|t_2-t_1|^{1/2}+|x_2-x_1|+|y_2-y_1|+|z_2-z_1|).$
The terminal condition $\phi$ is Lipschitz continuous.
\end{itemize}

\subsection{Dynamic programming equation}
We first consider a time discretization of the equation \eqref{edsr:general}. We denote the time step by $h=\frac{T}{N}$ and $(t_k=kh)_{0 \leq k \leq N}$ stand for the discretization times. For an arbitrary process $U$, set $\Delta U_k=U_{t_{k+1}}-U_{t_k}$. One needs to approximate $X$ by a Markov chain $X^N$ which can be simulated (take e.g. the Euler scheme \cite{jaco:04}). Whatever the scheme we use, we only require that it converges to $X$ in $L_2$:
\begin{itemize}
\item [\A{H3}] $\sup_{0\leq k\leq N} \e|X_{t_k}-\xn_{t_k}|^2 \rightarrow 0$ as $N$ goes to infinity.
\end{itemize}
The discrete-time counterpart of \eqref{edsr:general} is given by
$\yn_{t_N}=\phi(\xn_{t_N})$ and
 \begin{equation}
   \label{eq:yn}
   \yn_{t_k}  = \e_{t_k}(\yn_{t_{k+1}}) + h \e_{t_k} f(t_k,\xn_{t_k},\yn_{t_{k+1}},\zn_{t_k}) ,\qquad h\ \zn_{t_k}   =  \e_{t_k}(\yn_{t_{k+1}}\dw_{k}^*),
 \end{equation}
where $\e_{t_k}$ stands for the conditional expectation with
respect to $\f_{t_k}$ and $^*$ for the transpose. This is obtained by minimizing the difference $\e(\yn_{t_{k+1}}+h \e_{t_k} f(t_k,\xn_{t_k},\yn_{t_{k+1}},Z)-Y-Z\dw_k)^2$ over $\f_{t_k}$-measurable squared integrable random variables $(Y,Z)$. Our first result is the convergence of $(Y^N,Z^N)$ towards $(Y,Z)$ in the BSDE-norm, as $N$ goes to $\infty$. From $\yn$ and $\zn$, one could also easily define a process $L^N$ and prove its convergence to $L$. The result below seems to be new in this general setting.
\begin{theorem} \label{err:temps:discret}
Under \A{H1-H2-H3}, define the error
\begin{align*} e(N)=\max_{0 \leq k \leq N} \e |\yn_{t_k}-Y_{t_k}|^2 + \e \sum_{k=0}^{N-1}{\int_{t_k}^{t_{k+1}}{|\zn_{t_k}-Z_{t}|^2\dd t}},\end{align*}
where $\yn$ and $\zn$ are given by \eqref{eq:yn}. Then, $e(N)$ converges to $0$ as $N\rightarrow
\infty$. Furthermore, in the case of Brownian filtration
($\beta\equiv 0$ and $L\equiv 0$) and when $\xn$ is the Euler
scheme of $X$, one has $e(N)=O(N^{-1})$.
\end{theorem}
{\sc Proof.} In the sequel, $C$ denotes a constant which value changes throughout computations, but remains independent on $N$.

We begin by proving the result for $\max_{0 \leq k \leq N} \e |Y_{t_k}-Y^N_{t_k}|^2$. Firstly, we know \cite{elka:huan:97} that the solution
$(Y,Z,L)$ satisfies
\begin{equation}
  \label{eq:est:apriori}
  \e\big(\sup_{t\in[0,T]}Y^2_t+\int_0^T |Z_t|^2 \dd t +[L]_T\big)<+\infty.
\end{equation}
Then, from (\ref{edsr:general}-\ref{eq:yn}) we get
$Y_{t_k}-\yn_{t_k}=\e_{t_k}(Y_{t_{k+1}}-\yn_{t_{k+1}}) +
\e_{t_k}\int_{t_k}^{t_{k+1}}\{f(s,X_s,Y_s,Z_s) -\break
f(t_k,\xn_{t_k},\yn_{t_{k+1}},\zn_{t_k})\}\dd s.$ A combination of
Young's inequality $(a+b)^2\leq (1+\gamma h)a^2+(1+\frac1{\gamma h})b^2$ (with a parameter $\gamma>0$ to be chosen
later) and of the Lipschitz property of $f$ gives
\begin{align}
\nn \e|Y_{t_k}-\yn_{t_k}|^2  \leq & (1+\gamma h) \e|\e_{t_k}(Y_{t_{k+1}}-\yn_{t_{k+1}})|^2+C(h+\frac{1}{\gamma})\e \int_{t_k}^{t_{k+1}}{|Z_s-\zn_{t_k}|^2 \dd s} \\
 &+ C(h+\frac{1}{\gamma})(h^2+\int_{t_k}^{t_{k+1}} \e |X_s-\xn_{t_k}|^2 \dd s  +\int_{t_k}^{t_{k+1}} \e |Y_s-\yn_{t_{k+1}}|^2 \dd s). \label{app:eq1} \end{align}
Now define $\overline{Z}_{t_k}$ by $ h\overline{Z}_{t_k}:=\e_{t_k}\int_{t_k}^{t_{k+1}}{Z_s \dd s}=\e_{t_k}(\{Y_{t_{k+1}}+ \int_{t_k}^{t_{k+1}} f(s,X_s,Y_s,Z_s) \dd s\}\Delta W_k^*).$ Clearly
\begin{align}
\e \int_{t_k}^{t_{k+1}}{|Z_s-\zn_{t_k}|^2 \dd s} = \e
\int_{t_k}^{t_{k+1}}{|Z_s-\overline{Z}_{t_k}|^2 \dd s}
+h\e|\overline{Z}_{t_k}-\zn_{t_k}|^2. \label{app:eq2} \end{align}
The Cauchy-Schwarz inequality yields $|\e_{t_k}(\{Y_{t_{k+1}}-\yn_{t_{k+1}}\}\dw_{l,k})|^2 \leq h\{\e_{t_k}(|Y_{t_{k+1}}-\yn_{t_{k+1}}|^2)-|\e_{t_k}(Y_{t_{k+1}}-\yn_{t_{k+1}})|^2\}$ and consequently
\begin{align} \nn h\e|\overline{Z}_{t_k}-\zn_{t_k}|^2 \leq & C\e\{\e_{t_k}(|Y_{t_{k+1}}-\yn_{t_{k+1}}|^2)-\e|\e_{t_k}(Y_{t_{k+1}}-\yn_{t_{k+1}})|^2\} \\
& + C h\e \int_{t_k}^{t_{k+1}}{f(s,X_s,Y_s,Z_s)^2 \dd s} .
\label{app:eq3} \end{align} Plugging (\ref{app:eq2}-\ref{app:eq3})
into \eqref{app:eq1}, we get:
\begin{align*}
 \e|Y_{t_k}-\yn_{t_k}|^2
 \leq&  (1+\gamma h) \e|\e_{t_k}(Y_{t_{k+1}}-\yn_{t_{k+1}})|^2 +C(h+\frac{1}{\gamma})\e \int_{t_k}^{t_{k+1}}{|Z_s-\overline{Z}_{t_k}|^2 \dd s} \\
& +C(h+\frac{1}{\gamma})(h^2+\int_{t_k}^{t_{k+1}} \e |X_s-\xn_{t_k}|^2 \dd s  +\int_{t_k}^{t_{k+1}} \e |Y_s-\yn_{t_{k+1}}|^2 \dd s)\\
&  + C(h+\frac{1}{\gamma})\e\{\e_{t_k}(|Y_{t_{k+1}}-\yn_{t_{k+1}}|^2)-|\e_{t_k}(Y_{t_{k+1}}-\yn_{t_{k+1}})|^2\} \\
& + Ch(h+\frac{1}{\gamma})\e
\int_{t_k}^{t_{k+1}}{f(s,X_s,Y_s,Z_s)^2 \dd s}. \end{align*} Now
write $\e |Y_s-\yn_{t_{k+1}}|^2 \leq 2 \e |Y_s-Y_{t_{k+1}}|^2  + 2
\e |Y_{t_{k+1}}-\yn_{t_{k+1}}|^2 $ and analogously for
$X_{s}-\xn_{t_k}$, take $\gamma=C$: for $h$ small enough, it gives
\begin{align*}
 \e|Y_{t_k}-\yn_{t_k}|^2
 \leq & (1+C h) \e|Y_{t_{k+1}}-\yn_{t_{k+1}}|^2 + Ch^2+Ch\max_{0\leq k\leq N}\e |X_{t_k}-\xn_{t_k}|^2\\
&+C\e \int_{t_k}^{t_{k+1}}{|Z_s-\overline{Z}_{t_k}|^2 \dd s} +C\int_{t_k}^{t_{k+1}} \e |X_s-X_{t_k}|^2 \dd s\\
& +C\int_{t_k}^{t_{k+1}} \e |Y_s-Y_{t_{k+1}}|^2 \dd s + Ch\e
\int_{t_k}^{t_{k+1}}{f(s,X_s,Y_s,Z_s)^2 \dd s} \end{align*} and by
Gronwall's lemma $\max_{0\leq k\leq N} \e|Y_{t_k}-\yn_{t_k}|^2
 \leq  Ch+C\max_{0\leq k\leq N}\e |X_{t_k}-\xn_{t_k}|^2+C\sum_{k=0}^{N-1} \e
\int_{t_k}^{t_{k+1}}\{|Z_s-\overline{Z}_{t_k}|^2 + |X_s-X_{t_k}|^2
+ |Y_s-Y_{t_{k+1}}|^2\}\dd s.$
The contribution $\sum_{k=0}^{N-1}
\int_{t_k}^{t_{k+1}}\e|Y_s-Y_{t_{k+1}}|^2 \dd s$ is a $O(h)$:
indeed it is upper bounded by $3\sum_{k=0}^{N-1}
\int_{t_k}^{t_{k+1}} \dd s (t_{k+1}-s)\int_{s}^{t_{k+1}} \e
f(u,X_u,Y_u,Z_u)^2 \dd u + 3\sum_{k=0}^{N-1} \int_{t_k}^{t_{k+1}}
\dd s \int_{s}^{t_{k+1}} \e |Z_u|^2  \dd u + 3\sum_{k=0}^{N-1}
\int_{t_k}^{t_{k+1}} \dd s \e([L]_{t_{k+1}}-[L]_{s})$, which
equals a $O(h)$ owing to the a priori estimates \eqref{eq:est:apriori} on $(Z,L)$. In the same way, the contribution related
to $X_s-X_{t_k}$ is of order $O(h)$. Finally, it gives
\begin{align*}
\max_{0\leq k\leq N} \e|Y_{t_k}-\yn_{t_k}|^2 \leq &\ C h+C\max_{0\leq k\leq N}\e |X_{t_k}-\xn_{t_k}|^2+C\sum_{k=0}^{N-1} \e
\int_{t_k}^{t_{k+1}}|Z_s-\overline{Z}_{t_k}|^2 \dd s.
\end{align*}
Without extra assumptions, the above approximation's error related to the predictable process $Z$ converges to 0: combining this with \A{H3}, we conclude $\max_{0\leq k\leq N} \e|Y_{t_k}-\yn_{t_k}|^2\rightarrow 0$.\\
In the Brownian filtration case ($\beta\equiv 0$) and when $\xn$
is the Euler scheme of $X$, clearly $\e|X_{t_k}-\xn_{t_k}|^2=O(h)$ uniformly in $k$. Furthermore, Zhang
\cite{zhan:04} establishes that the error on $Z$ equals $O(h)$. Hence,
$\max_{0\leq k\leq N} \e|Y_{t_k}-\yn_{t_k}|^2=O(h)$.

Using the same techniques, it is easy to derive an estimate for $\e \sum_{k=0}^{N-1}{\int_{t_k}^{t_{k+1}}{|\zn_{t_k}-Z_{t}|^2\dd t}}$ from that on $\yn-Y$ (see also the proof of Theorem 2 in \cite{gobe:lemo:wari:05}). We omit further details.\qed

\subsection{Description of the algorithm}
\label{par:algo}
At some point in our approach, the numerical solution of $(Y,Z)$ needs to be upper bounded, especially to analyze the empirical regression errors. This is a theoretical reason but we also mention that it seems to have some importance in our numerical experiments. To get a bounded solution, we consider a threshold $R=(R_0,R_1,\ldots,R_d)\in (\r^+)^{d+1}$ (a priori with large coordinates) and this threshold defines
\begin{align} [\dw_{l,k}]_w&=\big(-R_0\sqrt{h}\ \big) \vee \dw_{l,k}\wedge \big(R_0\sqrt{h}\ \big),\nn\\
\nonumber f^R(t,x,y,z) & = f(t,-R_1 \vee x_1\wedge R_1,\cdots,-R_{d} \vee x_{d} \wedge R_{d},  y,z), \\
 \nn \phi^{R}(x) & = \phi(-R_1 \vee x_1\wedge R_1,\cdots,-R_{d} \vee x_{d} \wedge R_{d}),\\
\nn [\psi]_y(x)&=-C_y(R) \vee \psi(x) \wedge C_y(R),\quad [\psi]_z(x)=-\frac{C_y(R)}{\sqrt h} \vee \psi(x) \wedge \frac{C_y(R)}{\sqrt h}.
\end{align}
The constant $C_y(R)$ is explicit and defined later (see equation \eqref{eq:cy} below).

It is easy to see that the conditional expectations in \eqref{eq:yn} are regression functions of $X^N_{t_k}$; hence, it is natural to approximate them by their projections on finite-dimensional bases. For the scalar component of $Y$ and the $q$ components of $Z$, at each discretization time $t_k$, we take $q+1$ deterministic function bases $\big(p_{l,k}(\cdot)\big)_{0 \leq l \leq q}$. For sake of convenience, the basis $p_{l,k}$ of size $K_{l,k}$ is considered as a $K_{l,k}$-dimensional vector of scalar functions. A function $F$ in the vector space spanned by the basis $p_{l,k}$ is given by its coefficients $\alpha$ and one uses the short notation $F(\cdot)=\alpha\cdot p_{l,k}(\cdot)$.

We denote by $(\xnm_{t_k})_{1 \leq m \leq M,0 \leq k \leq N}$ and $(\dw_k^m)_{1 \leq m \leq M,0 \leq k \leq N-1}$ the $M$ independent simulations of $(\xn_{t_k})_{0 \leq k \leq N}$ and $(\dw_k)_{0 \leq k \leq N-1}$. \\

With these notations, we are in a position to define $\pynrm_k(\xnk)$ and $\pznrm_{l,k}(\xnk)$ (the approximations of $Y^N_{t_k}$ and $Z^N_{l,t_k}$), where the functions $\pynrm_k(\cdot)$ and $\pznrm_{l,k}(\cdot)$ are given as follows.
\begin{itemize}
\item [$\pmb \rightarrow$] Initialization : for $k=N$ set $\pynrm_N(\cdot)=\phi^{R}(\cdot)$.
\item [$\pmb \rightarrow$] Iteration : for $k=N-1, \cdots, 0$, solve the $q$ least-squares problems\footnote{For the numerical resolution of least-squares problems, see \cite{golu:vanl:96}.}:
\[\alpha_{l,k}^M = \arg\inf_{\alpha}
\frac{1}{M}\sum_{m=1}^{M}{|\pynrm_{k+1}(\xnm_{t_{k+1}})\frac{[\dw_{l,k}^m]_w}{h}-\alpha\cdot p_{l,k}(\xnm_{t_{k}})|^2} . \]
Put \fbox{$\pznrm_{l,k}(\cdot)=[\alpha_{l,k}^M\cdot p_{l,k}(\cdot)]_z$}. Then compute $\alpha_{0,k}^M$ as the minimizer of
\[\hspace{-1cm} \frac{1}{M}\sum_{m=1}^{M}|\pynrm_{k+1}(\xnm_{t_{k+1}})+h f^R(t_k,\xnm_{t_k},\pynrm_{k+1}(\xnm_{t_{k+1}}), \pznrm_{l,k}(\xnm_{t_{k}})) - \alpha \cdot p_{0,k}(\xnm_{t_{k}})|^2 \]
over $\alpha$. In the above definition, we write
$f^R(t_k,x,y,z_l)$ for $f^R(t_k,x,y,(z_l)_{1\leq l\leq q})$.
Put \fbox{$\pynrm_{k}(\cdot)=[\alpha_{0,k}^M\cdot p_{0,k}(\cdot)]_y$}.
\end{itemize}

\subsection{Convergence analysis}
\subsubsection{Error from the threshold $R$}
\label{par:erreur localisation}
We first compare the approximative solution $(Y^N,Z^N)$ defined by \eqref{eq:yn} with $(\ynr,\znr)$ defined by:
\begin{equation}
\label{eq:ynr} \ynr_{t_k}=  \e_{t_k}(\ynr_{t_{k+1}}) + h\e_{t_k}f^R(t_k,\xn_{t_k},\ynr_{t_{k+1}},\znr_{t_k}) , \qquad h\znr_{t_k} = 
\e_{t_k}(\ynr_{t_{k+1}}[\dw_k]^*_w), \end{equation} and
$\ynr_{t_N}=\phi^R(\xn_{t_N})$. In \cite{lemo:gobe:wari:05}, we prove that $(\ynr,\znr)$ converges in $L_2$ to $(\yn,\zn)$ as $R$ goes to infinity. A rate of convergence is also available under the additional assumption $\sup_{0\leq k\leq N}\e|\xn_{t_k}|^p\leq C_p(1+|x|^p)$ for some $p>2$. Namely, provided that $R_i= N^{2/(p-2)}$ ($i=1,\cdots,d$) and $R_0=c\sqrt{\log(N)}$ (for $c$ large enough), one has
$\max_{0 \leq k \leq N} \e |\ynr_{t_k}-\yn_{t_k}|^2 + h\e\sum_{k=0}^{N-1}{|\znr_{t_k}-\zn_{t_k}|^2}=O(N^{-1})$ as in Theorem \ref{borne erreur}. Hence, especially if $p$ is not small, reasonable values of $R$ yield a very good approximation.

The new dynamic programming equation \eqref{eq:ynr} has the main advantage to give bounded solutions. Namely one has $|\ynr_{t_k}| \leq  C_y(R)$ and $\sqrt h |\znr_{t_k}|\leq C_y(R)$ where 
\begin{equation}
  \label{eq:cy}
C_y(R) =  e^{(2\gamma^*+\frac{1+\gamma^*}{q} )T} \left\{  \sup_{x}|\phi^{N,R}(x)|^2 +  2T\frac{1+\gamma^*}{\gamma^*} \sup_{t,x}|f^R(t,x,0,0)|^2 \right\}
\end{equation}
  and $\gamma^*=4qC_f^2$ (see \cite{lemo:gobe:wari:05} and \cite{lemo:05}). Note that the $L_\infty$-norms related to $\phi^{N,R}$ and $f^R$ can be easily computed on usual examples.

\subsubsection{Error from the projections on bases and the simulations}
\label{par:h4}
The following extra assumption \A{H4} is sufficient to prove that the functions $y^{N,R}_k(\cdot)$ and $\sqrt h z^{N,R}_k(\cdot)$ (defined by  $y_k^{N,R}(\xn_{t_k})=Y^{N,R}_{t_k}$ and $z_k^{N,R}(\xn_{t_k})=Z^{N,R}_{t_k}$) are Lipschitz continuous, uniformly in $h$ and $R$ (see \cite{lemo:gobe:wari:05}).
\begin{itemize}
\item [\A{H4}] Denote by $X^{N,k_0,x}_{t_k}$ the Markov chain $\xn_{t_k}$ starting at $x$ at time $t_{k_0}$. There is a constant $C>0$ such that
  \begin{itemize}
  \item [a)] $\e |X^{N,k_0,x}_{t_N} - X^{N,k_0,x'}_{t_N} |^2 + \e  |X^{N,k_0,x}_{t_{k_0+1}} - X^{N,k_0,x'}_{t_{k_0+1}} |^2\leq  C|x-x'|^2$ for any $x$ and $x'$, uniformly in $k_0$ and $N$.
\item [b)] $\e|X^{N,k_0,x}_{t_{k_0+1}} - x|^2\leq
  C h (1+|x|^2)$ for any $x$, uniformly in $k_0$ and $N$.

  \end{itemize}
\end{itemize}
This assumption on $X^N$ is not restrictive since this property holds for $X$ under \A{H1} (and usually, for the Euler scheme as well). The next result states a general upper bound for the error between $(\ynr_{t_k},\znr_{t_k})$ and the solutions $(\pynrm_k(\xnk),\pznrm_k(\xnk))$ computed with our scheme, in terms of the basis functions and the number of simulations. Specifications of the upper bounds are given later, according to the choice of bases.
\begin{theorem} \label{borne erreur} (see \cite{lemo:gobe:wari:05}) Assume \A{H1-H2-H3-H4} and denote by $K^M_{l,k}$ the random variable given by the rank of the matrix of size $K_{l,k}\times M$ which columns are $(p_{l,k}(X^{N,m}_{t_k}))_m$ (note that $K^M_{l,k}\leq K_{l,k}$). By convention, we set $K_{0,N}=0$. \\ Then, there exists a constant $C$ such that for any $\beta\in]0,1]$ one has
\begin{align*} & \max_{0 \leq k \leq N} \e |\pynr_k(\xn_{t_k})-\pynrm_k(\xn_{t_k})|^2   + h \e \sum_{k=0}^{N-1}|\pznr_k(\xn_{t_k})-\pznrm_k(\xn_{t_k})|^2 \\
\leq & C\frac{C_y(R)^2\log(M)}{M}\sum_{k=0}^{N-1}{\sum_{l=0}^{q}K_{l,k}} + Ch^{\beta} \\
& + C\sum_{k=0}^{N-1}{\{\inf_{\alpha} \e |\pynr_k(\xn_{t_k})-\alpha \cdot p_{0,k}(\xn_{t_k})|^2  + \sum_{l=1}^{q}{  \inf_{\alpha} \e |\sqrt{h}\pznr_{l,k}(\xn_{t_k})-\alpha \cdot p_{l,k}(\xn_{t_k})|^2  }\}} \\
& + C \frac{C_y(R)^2}{h} \sum_{k=0}^{N-1}\bigg\{ \e \big( K_{0,k}^M \exp(-\frac{Mh^{\beta+2}}{72C_y(R)^2K_{0,k}^M}) \exp(CK_{0,k+1}\log\frac{C\ C_y(R)(K_{0,k}^M)^{\frac{1}{2}}}{h^{\frac{\beta+2}{2}}}) \big)\\
& + h\e \big( K_{l,k}^M \exp(-\frac{Mh^{\beta+1}}{72C_y(R)^2R_0^2 K_{l,k}^M})  \exp(CK_{0,k+1}\log\frac{C\ C_y(R)R_0(K_{l,k}^M)^{\frac{1}{2}}}{h^{\frac{\beta+1}{2}}}) \big)\\
&+ \exp(C K_{0,k}\log\frac{C\ C_y(R)}{h^{\frac{\beta+2}{2}}})
\exp(-\frac{Mh^{\beta+2}}{72C_y(R)^2})  \bigg\}. \end{align*}
\end{theorem}
The parameter $\beta$ should be chosen optimally to get the sharpest upper bound, see Paragraph \ref{par:complexity}. Theorem \ref{borne erreur} improves our previous results \cite{gobe:lemo:wari:05}, where the statistical errors are estimated in terms of the fourth moments of the $L_2$-orthonormalized basis functions. These quantities are generally unknown, which makes the adjustment of the bases, $h$ and $M$ difficult. In order to interpret terms of the upper bound above, we first recall a well-known result in non-parametric regressions.
\begin{theorem} \label{th:gyorfi} (see \cite{gyor:kohl:krzy:walk:02} Theorem 11.1 p. 184). Let $(U,V)$ be two random variables taking values in $\Re^d\times\Re$ and let $F$ be a $K$-dimensional linear space of functions. Consider the problem of estimating the regression function $v(u)=\e(V|U=u)$, using $M$ independent realizations $(U_m,V_m)_{1 \leq m \leq M}$ of $(U,V)$. Approximate $v(\cdot)$ by its best approximation in $F$ using an empirical $L_2$-norm related to the $M$ observations:
\[\hat{v}_M(\cdot) = \arg \inf_{f \in F} \frac{1}{M}\sum_{m=1}^{M}{|V_m-f(U_m)|^2} .\]
Under the assumption $\sigma^2=\sup_{u} {\rm Var}(V|U=u) < \infty$, one has
\begin{equation}
\label{eq:distribution:free}
 \e \bigg( \frac{1}{M}\sum_{m=1}^{M}{|v(U_m) - \hat{v}_M(U_m)|^2} \bigg) \leq  \sigma^2 \frac{K}{M} + \min_{f \in F}\e|f(U)-v(U)|^2 . \end{equation}
\end{theorem}
This is the prototype of {\em model-free} results since one only requires the conditional variance to be uniformly bounded. 

In our case, at each $t_k$ one solves a regression problem, which yields errors of type \eqref{eq:distribution:free} (compare with the first and second lines of the r.h.s. of Theorem \ref{borne erreur}). The boundedness of $(\ynr,\sqrt h \znr)$ ensures that the relative conditional variances are uniformly bounded. The additional factor $\log(M)$ in the first line occurs because one estimates $\e |\pynr_k(\xn_{t_k})-\pynrm_k(\xn_{t_k})|^2$ instead of $\e \frac1M \sum_{m=1}^M|\pynr_k(\xnm_{t_k})-\pynrm_k(\xnm_{t_k})|^2$ as in \eqref{eq:distribution:free}.

However, the scheme to currently analyze is much more delicate than in Theorem \ref{th:gyorfi}, because the regression problems are all correlated: indeed, $V_m$ should be given by $\pynrm_{k+1}(\xnm_{t_{k+1}})$ which is not independent of $V_{m'}$ ($m\neq m'$) because of the random function $\pynrm_{k+1}(\cdot)$. The trick consists in covering the class of functions where $\pynrm_{k+1}(\cdot)$ is, by a finite number of balls centered at deterministic functions. The radius of the above balls depends on $h$ and on an extra parameter $\beta$. The covering number is estimated using the Vapnik-Chervonenkis dimension. All these extra contributions are the other terms in the r.h.s. of Theorem \ref{borne erreur}.\\
The analysis of the algorithm's error turns out to be quite intricate and actually, we guess that our estimates of Theorem \ref{borne erreur} are not optimal. However, we do not have good ideas to improve them. Before performing an analysis of complexity in Paragraph \ref{par:complexity}, we give a modified algorithm which complexity is the same and which theoretical upper bound is better (paradoxically, both schemes behave similarly on numerical tests).

\subsubsection{A modified algorithm}
\label{par:algo:modif} As mentioned before, some terms in the upper bound of Theorem \ref{borne erreur} come from the lack of independence between least-squares problems at different discretization times $t_k$. These terms are those from the third and fourth lines in the r.h.s. of the upper bound. We present below a variant of the algorithm of Paragraph \ref{par:algo}, which removes these two terms (for details, see \cite{lemo:05}).

\begin{itemize}
\item [$\pmb \rightarrow$] Initialization : for $k=N$ set $\pynrm_N(\cdot)=\phi^{R}(\cdot)$.
\item [$\rightarrow$] Let $k < N-1$. For each path $m$, simulate $(\xntm_{t_{k+1}},\dwt_k^m)$ which are, conditionally to $\xnm_{t_k}$, an independent copy of $(\xnm_{t_{k+1}},\dw_k^m)$ (and independent of everything else). The coefficients $\alpha_{l,k}^M$ and $\alpha_{0,k}^M $ are now computed as 
\begin{align*}
\alpha_{l,k}^M &= \arg\inf_{\alpha}\frac{1}{M}\sum_{m=1}^{M}{|\pynrm_{k+1}(\xntm_{t_{k+1}})\frac{[\dwt_{l,k}^m]_w}{h}-\alpha
\cdot p_{l,k}(\xnm_{t_{k}})|^2},\\
\alpha_{0,k}^M& = \arg\inf_{\alpha} \frac{1}{M}\sum_{m=1}^{M}|\pynrm_{k+1}(\xntm_{t_{k+1}})\\
&+h f^R(t_k,\xnm_{t_k},\pynrm_{k+1}(\xntm_{t_{k+1}}), \pznrm_{l,k}(\xnm_{t_{k}}))  - \alpha \cdot p_{0,k}(\xnm_{t_{k}})|^2.
\end{align*}
As before, the functions $\pynrm_{k}(\cdot)$ and $\pznrm_{l,k}(\cdot)$ are given by $[\alpha_{0,k}^M\cdot  p_{0,k}(\cdot)]_y$ and $[\alpha_{l,k}^M\cdot  p_{l,k}(\cdot)]_z$.
\end{itemize}
Note that this algorithm requires to draw only twice more simulations.

\subsection{Accuracy and complexity}\label{par:complexity}
Now we focus on the error between the discrete-time dynamic programming equation \eqref{eq:yn} and the solutions computed by our initial algorithm, its modified version and also the Bouchard-Touzi's algorithm \cite{bouc:touz:04}. Our aim is at discussing the trade-off between the accuracy and the complexity (i.e. computational time).
\subsubsection{Initial algorithm (see Paragraph \ref{par:algo})}
\label{par:HC}
From Theorem \ref{borne erreur}, we can derive how to make
$N$, $M$ and the number of basis functions vary together. In this
discussion, we neglect the influence of the localization parameter
$R$, which is supposed to be large enough from the beginning (see
paragraph \ref{par:erreur localisation}). As
already observed in \cite{gobe:lemo:wari:05}, local basis functions
take advantage of the Lipschitz property of functions
$\pynr_k(\cdot)$ and $\sqrt{h}\pznr_{l,k}(\cdot)$ (see Paragraph \ref{par:h4}). Let us consider the simplest
example of a local function basis, i.e. the hypercubes basis,
already used in \cite{gobe:lemo:wari:05} and still denoted {\bf HC}
here. To simplify, $p_{l,k}$ does not depend on $l$ or $k$ and its
size equals $K$. Choose a domain $D\subset \r^{d}$ centered
on $x$, that is $D=\prod_{i=1}^{d}{]x_i-a,x_i+a]}$, and partition
it into small hypercubes of edge $\delta$. Thus,
$D=\cup_{i_1,\cdots,i_{d}} D_{i_1,\cdots,i_{d}}$ where
$D_{i_1,\cdots,i_{d}}=]x_1-a+i_1\delta,x_1-a+(i_1+1)\delta] \times
\cdots \times ]x_d-a+i_{d}\delta,x_d-a+(i_{d}+1)\delta]$. Then we
define $p_{l,k}$ as the indicator functions associated to this set
of hypercubes:
$p_{l,k}(\cdot)=\big(\1_{D_{i_1,\cdots,i_{d}}}(\cdot)
\big)_{i_1,\cdots,i_{d}}$. With this particular choice of function
bases, we can make the projection error of Theorem \ref{borne
erreur} explicit and refer to \cite{gobe:lemo:wari:05} for details:
\begin{align*}
\inf_{\alpha}\e \big(|\pynr_k(\xn_{t_k})-\alpha \cdot
p_{0,k}(\xn_{t_k})|^2 \big) \leq \ C \{\delta^2 + C_y(R)^2 \p(
\xn_{t_k} \in D^c) \} .\end{align*}
As for the impact of the threshold $R$, $\p( \xn_{t_k} \in D^c)$ becomes negligible with respect to the other errors if we choose $D$ big enough (a feature which is confirmed by the numerical experiments). Thus and as far the projection errors are concerned, to get a global (squared) error of order $h^{\beta}$ we have to choose $\delta \approx h^{\frac{\beta+1}{2}}$, or equivalently a number of basis functions $K \approx h^{-\frac{d(\beta+1)}{2}}$ (considering a fixed domain $D$). Regarding now the number of simulations $M$, to avoid an explosive upper bound in Theorem \ref{borne erreur}, one should take $M \approx Ch^{-d(\beta+1)-(\beta+2)}\log(h^{-\frac{d(\beta+1)}{4}-\frac{\beta+1}{2}})$ for a constant $C$ large enough (here, the ranks $K_{l,k}^M$ have simply been upper bounded by $K$).\\
The dominant term of the algorithm's complexity $\mc{C}$ associated to this choice of function basis is $\mc{C}=NMd\log(K)$, which corresponds to determine in which cells the simulation fall (this is the cost of a nearest neighbor algorithm in a tensored grid, i.e. $O( d \log(K))$ for one path at a given time). Hence, up to logarithmic factors, the complexity equals $\mc{C}=O(h^{-1-d(\beta+1)-(\beta+2)})$, while the squared error is of order $h^\beta=O(\mc{C}^{-\beta/(2+(\beta+1)(d+1))})$. The optimal value of $\beta\in]0,1]$ is $\beta=1$, for which the squared error is of order $h=O(\mc{C}^{-1/(2 d +4)})$, while $\delta\approx h$ and $M\approx Ch^{-2d-3}\log(1/h)$.
\subsubsection{Modified algorithm (see Paragraph \ref{par:algo:modif})}
The complexity is unchanged. But the error bounds are better since the terms from the third and fourth lines in the r.h.s. of Theorem \ref{borne erreur} have disappeared. Hence the asymptotics on $M$ is less stringent: indeed, we get $M\approx Ch^{-d-3}\log(1/h)$, $\delta\approx h$ in order to achieve a squared error of order $h=O(\mc{C}^{-1/(d +4)})$.
\subsubsection{Bouchard-Touzi's algorithm \cite{bouc:touz:04}}
We compute the complexity of their algorithm in the most favorable case where $X$ is a geometric Brownian motion. Otherwise in a more general diffusion framework, the algorithm is heavy to implement and its complexity more difficult to evaluate because of the necessary calculation of Skorohod's integrals. \\
In this algorithm, one needs $N$ independent sets $(\mc{M}_k)_{1 \leq k \leq N}$ of  $M$ simulated paths of $X^N$ (one set for each discretization time). At each discretization time $t_k$ and for each path in the set $\mc{M}_k$, a calculation involving the $M$ paths of the set $\mc{M}_{k+1}$ is performed. This leads to a complexity $\mc{C}=O(NM^2)$. The squared error associated to this complexity is given by Theorem 6.2 in \cite{bouc:touz:04} and is of order $\frac 1N+\frac{N^{2+\frac{d}{4}}}{\sqrt{M}}$. Expressing the squared error as a function of the complexity, we find $\mc{C}^{-\frac{1}{13+d}}$.
\subsubsection{Summary}
\begin{table}[h]
\begin{center}
\begin{tabular}{|c|c|c|}
\hline
Initial Algorithm & Modified Algorithm & Bouchard-Touzi's Algorithm \\
(Paragraph \ref{par:algo})& (Paragraph \ref{par:algo:modif})& \cite{bouc:touz:04}\\
\hline
   &   & $\mc{C}^{-\frac{1}{13+d}}$\\$\mc{C}^{-\frac{1}{4+2 d}}$ &$\mc{C}^{-\frac{1}{4+d}}$& (if $X$=Brownian motion\\
&& or geometric BM)\\
\hline
\end{tabular}
\caption{Squared error with respect to
the complexity $\mc{C}$.\label{table:complexity}}
\end{center}
\end{table}
It turns out that in the geometric Brownian case, our algorithm is more efficient than Bouchard-Touzi's Algorithm for $d \leq 9$ and less efficient otherwise. But the geometric Brownian framework is very favorable for a Malliavin calculus approach, since the integration by parts formulas are especially simple and this is no more true for general diffusion models. With our approach, the complexity is independent of the model.

Furthermore, the modified algorithm is very competitive in any case, although it is less natural (because of the extra simulations). We have observed in practice that the initial and modified algorithm differ very little from each other. This suggests that the upper bound in Theorem \ref{borne erreur} may be not optimal (and the resulting complexity as well). We hope that we could address these issues in future works.

\section{Reflected backward stochastic differential equation}
 \label{section:2}

In this section, we restrict to a Brownian BSDE reflected on one barrier. We show that the results presented in the non-reflected case can apply to the reflected case. Actually we present three approximation methods, which we call in short {\em max method}, {\em penalization method} and {\em regularization method}.

The framework is analogous to Section \ref{section:1}, without the Poisson random measure process. Namely, we suppose that the forward process is a diffusion SDE $X_{t}=x+\int_{0}^{t}{b(s,X_s) \dd s} + \int_{0}^{t}\sigma(s,X_s)\dd W_s$ under the following assumption
\begin{itemize}
\item [\A{H1'}] The functions $b$ and $\sigma$ are continuously differentiable with uniformly bounded derivatives w.r.t. $x$. The matrix $\sigma\sigma^*$ satisfies the ellipticity condition $\sigma \sigma^* \geq \epsilon \id$ with $\epsilon>0$.
\end{itemize}
Regarding the reflected BSDE (RBSDE in short), we consider a continuous function $\phi$ which defines a continuous adapted obstacle $\phi(t,X_t)$. The RSBDE $(Y,Z,K)$ is solution of
\begin{equation}
  \label{rbsde}
  \left\{\begin{array}{l}
Y_t =  \phi(T,X_T)+\int_t^T f(s,X_s,Y_s,Z_s)\dd s+{K_T-K_t}-\int_t^T Z_s \dd W_s,\\
Y_t \geq \phi(t,X_t) ,\\
\text{$K$ is continuous, increasing, $K_0=0$ and }\int_0^T (Y_t-\phi(t,X_t) )\dd K_t=0.   
\end{array}\right.
\end{equation}
To ensure existence/uniqueness of the process above \cite{elk:97} and further approximation results, we impose stronger conditions on the driver and the obstacle:
\begin{itemize}\item [\A{H2'}] For any $(t_1,x_1,y_1,z_1), (t_2,x_2,y_2,z_2)\in [0,T]\times \r^d\times \r \times \r^q$, one has $|f(t_2,x_2,y_2,z_2)-f(t_1,x_1,y_1,z_1)|\leq C_f(|t_2-t_1|^{1/2}+|x_2-x_1|+|y_2-y_1|+|z_2-z_1|)$
 and $|\phi(t_2,x_2)-\phi(t_1,x_1)|\leq C_\phi(|t_2-t_1|^{1/2}+|x_2-x_1|)$.
\end{itemize}
When the driver $f$ is linear, i.e. $f(t,x,y,z)=-r_t y -z \ \theta_t$ where $\theta$ is the market risk-premium in finance, it is well known that $Y_t$ equals the price of the American option with payoff $P_t=\phi(t,X_t)$: $Y_t={\rm ess\ sup}_{\text{stopping times }\tau \in[t,T]}\e_{\mathbb{Q}}(e^{-\int_t^\tau r_s \dd s} P_\tau|{\cal F}_t).$
\subsection{Approximation procedures}

We denote by $\xn$ the Euler scheme of $X$, based on the time net $(t_k=kh)_{0 \leq k \leq N}$ where $h=\frac{T}{N}$. As before, one needs $M$ independent simulations of $(\xn,\Delta W)$. We propose three procedures to approximate $(Y,Z)$, solution of the reflected BSDE above.

\subsubsection{Max method}
This method is standard and simply consists in taking at each time $t_k$ the maximum between the obstacle and the expected value of the non reflected BSDE. The resulting approximation $(\yn,\zn)$ is defined by:
\begin{align*}
\begin{cases}
\yn_{t_N} = & \phi(t_N,\xn_{t_N}) ,\\
\yn_{t_k} = & \max\big(\phi(t_k,\xn_{t_k}),\e_{t_k}( \yn_{t_{k+1}} + hf(t_k,\xn_{t_k},\yn_{t_{k+1}},\zn_{t_k}))\big) , 0 \leq k \leq N-1, \\
  \zn_{t_k} = & \frac 1 h\e_{t_k}(\yn_{t_{k+1}} \dw^*_{k}), 0 \leq k \leq N-1. \end{cases}
\end{align*}
Ma and Zhang \cite{ma:zhan:05} prove that the rate of convergence of such approximation is $N^{-1/4}$ when the obstacle function is smooth (i.e. of class $C^{1,2}$). In \cite{lemo:05} we extend this analysis to less smooth obstacles in order to handle for instance the case of American put contract (i.e. $\phi(t,x)=(K-e^x)_+$ if $X$ is the log-price process). The relevant assumption of $\phi$ is the following one:
\begin{itemize}\item [\A{H3'}] The function $\phi$ is of class $C^1$ w.r.t. the time variable and uniformly Lipschitz continuous  w.r.t. the space variable. Moreover, $\phi(t,X_t)$ and $\phi(t,X^N_t)$ satisfy the following It\^o expansions:
\begin{align*} \phi(t,X_t)=& \phi(0,x)+\int_{0}^{t}{U_s \dd s} + \int_{0}^{t}{V_s \dd W_s} + A_t, \\
\phi(t,\xn_t) = & \phi(0,x)+\int_{0}^{t}{U^N_s \dd s} + \int_{0}^{t}{V^N_s \dd W_s} + A^N_t, \end{align*}
where $A$, $A^N$ are increasing, continuous, integrable and such that the measures $\dd A_t$ and $\dd A^N_t$ are singular  w.r.t. $\dd t$ and where $U$, $U^N$, $V$ and $V^N$ satisfy the following integrability condition $\sup_{t\le T} \e(|U_t|^p)+\sup_{t\le T,N} \e ( |U^N_t|^p ) + \sup_{t\le T} \e (|V_t|^p) + \sup_{t\le T,N} \e (|V^N_t|^p ) < \infty$ for $p \geq 2$.
\end{itemize}
\begin{theorem} (see \cite{lemo:05}) Under Assumption \A{H1'-H2'-H3'}, one has
\begin{align*} \max_{0 \leq k \leq N} \e  |\yn_{t_k}-Y_{t_k}|^2 + \sum_{k=0}^{N-1}{\int_{t_k}^{t_{k+1}}{\e  |\zn_{t_k}-Z_t|^2  dt}} \leq C(1+|x|^4)N^{-1/2}. \end{align*}
\end{theorem}
Note that the convergence rate is lower than in the non-reflected case ($N^{-\frac{1}{4}}$ instead of $N^{-\frac{1}{2}}$). The algorithm to approximate $(\yn,\zn)$ is very similar to the one of Paragraph \ref{par:algo} (with analogous threshold functions).
\begin{enumerate}
\item $y^{N,R,M}_N(\cdot)=\phi^R(t_N,\cdot)$
\item At a given discretization time $t_k$ ($0 \leq k \leq N-1$), the coefficients $\alpha_{l,k}^M$ ($1\leq l \leq q$) are defined as the minimizers over $\alpha_l$ of $\displaystyle\frac{1}{M}\sum_{m=1}^{M}{|y_{k+1}^{N,R,M}(\xnm_{t_{k+1}})\frac{[\dw_{l,k}^m]_w}{h}-\alpha_l.p_{l,k}^m|^2}.$
Then put $z_{l,k}^{N,R,M}(\cdot)=[\alpha_{l,k}^M.p_{l,k}]_z(\cdot)$.
\item $\alpha_{0,k}^M$ is finally defined as the minimizer over $\alpha_0$ of
\begin{align*}
\frac{1}{M}\sum_{m=1}^{M}{|y^{N,R,M}_{k+1}(\xnm_{t_{k+1}})+hf^R(t_k,\xnm_{t_k},y_{k+1}^{N,R,M}(\xnm_{t_{k+1}}),z^{N,R,M}_{k}(\xnm_{t_k}))-\alpha_0.p_{0,k}^m|^2}. \end{align*}
Then put $y^{N,R,M}_k(\cdot)=\max(\phi^R(t_k,\cdot),[\alpha_{0,k}^M.p_{0,k}]_y(\cdot))$.
\end{enumerate}
The error between $(\ynr_{t_k},\znr_{t_k})$ and $(y_k^{N,R,M}(\xn_{t_k}),z_k^{N,R,M}(\xn_{t_k}))$ can be analyzed exactly as in Theorem \ref{borne erreur} and actually, we get the same estimates.

\subsubsection{Penalization method}
This method is based on the ideas developed in \cite{elk:97} to show the existence and uniqueness of the solution to \eqref{rbsde}: the reflection is handled via a penalization term which is very large only when $Y$ is below to the obstacle. Namely, the solution $(Y,Z,K)$ is approximated by $(Y^n,Z^n,K^n)$ where $(Y^n,Z^n)$ is the solution of the standard BSDE 
\begin{align*} 
Y^n_t = \phi(T,X_T) + \int_{t}^{T}{f(s,X_s,Y^n_s,Z^n_s)\dd s} + n \int_{t}^{T}{(Y^n_s - \phi(s,X_s))_-\dd s} - \int_{t}^{T}{Z^n_s \dd W_s}, \end{align*}
and $K^n$ is defined by $K^n_t=n\int_{0}^{t}{(Y^n_s-\phi(s,X_s))_- \dd s}$. It is thus shown in \cite{elk:97} that $(Y^n,Z^n,K^n)$ tends to $(Y,Z,K)$ as $n$ tends to infinity with the property $Y^n \uparrow Y$. In general, one does not know the convergence rate (some results are available when $f$ does not depend on $Z$, see \cite{hama:jean:05}). At a fixed $n$, $(Y^n,Z^n)$ is the solution of a standard BSDE and can thus be approximated by the algorithm proposed in the first section:
\begin{enumerate}
\item $y^{n,N,R,M}_N(\cdot)=\phi^R(t_N,\cdot)$
\item At a given discretization time $t_k$ ($0 \leq k \leq N-1$), the coefficients $(\alpha_{l,k}^M)$ are defined as the minimizers of $\displaystyle \inf_{\alpha_l} \frac{1}{M}\sum_{m=1}^{M}{|y_{k+1}^{n,N,R,M}(\xnm_{t_{k+1}})\frac{[\dw_{l,k}^m]_w}{h}-\alpha_l.p_{l,k}^m|^2}.$
Then put $z_{l,k}^{n,N,R,M}(\cdot)=[\alpha_{l,k}^M.p_{l,k}]_z(\cdot)$. 
\item $\alpha_{0,k}^M$ is finally defined as the minimizer of:
\begin{align*}
\inf_{\alpha_0} \frac{1}{M}\sum_{m=1}^{M}|& y^{n,N,R,M}_{k+1}(\xnm_{t_{k+1}})+hf^R(t_k,\xnm_{t_k},y_{k+1}^{n,N,R,M}(\xnm_{t_{k+1}}),z^{n,N,R,M}_{l,k}(\xnm_{t_k})) \\
& +nh(y^{n,N,R,M}_{k+1}(\xnm_{t_{k+1}})-\phi^R(t_k,\xnm_{t_k}))_- -\alpha_0.p_{0,k}^m|^2.\end{align*}
Then put $y^{n,N,R,M}_k(\cdot)=[\alpha_{0,k}^M.p_{0,k}]_y(\cdot)$.
\end{enumerate}

\subsubsection{Regularization method}
The last idea is based on the results of \cite{ball:02}, which exploit the specific form of the increasing process $K$. The authors consider the case where the obstacle $\phi(t,X_t)$ is a c\`adl\`ag semimartingale satisfying $\phi(t,X_t)=\phi(0,x)+\int_{0}^{t}{U_s \dd s} + \int_{0}^{t}{V_s \dd W_s} + A_t$ 
where $\e \int_{0}^{T}{\{|V_t|^2+|U_t|^2\}dt} < \infty$ and where $A$ is a c\`adl\`ag adapted increasing process with $\e |A_T| < \infty$. In this case, a solution $(Y,Z,K)$ to \eqref{rbsde} is also a solution to
\begin{align*} \begin{cases} Y_t= & \phi(t,X_t) + \int_{t}^{T}{\{f(s,X_s,Y_s,Z_s) + \alpha_s \1_{Y_s=\phi(s,X_s)} [ f(s,X_s,\phi(s,X_s),V_s)+U_s]_- \}\dd s} \\
& - \int_{t}^{T}{Z_s \dd W_s}, \\
Y_t \geq & \phi(t,X_t) \end{cases} \end{align*}
with $\e \int_{0}^{T}{ \{ |Y_t|^2 + |Z_t|^2 + |\alpha_t|^2 \} } dt < \infty$.
To show the existence and uniqueness of the solution to the equation above, the authors use a regularization method based on $C^{\infty}$ functions $\varphi_n$ (with $0 \leq \varphi_n \leq 1$) satisfying $\varphi_n(x)=1$ if $|x| \leq \frac{1}{2^n}$ and $\varphi_n(x)=0$ if $|x| \geq \frac{2}{2^n}$, which are mollifiers of $\varphi(x)=\1_{x=0}$. It turns out that the solution $(Y^n,Z^n)$ to the following BSDE converges to $(Y,Z)$ as $n$ goes to infinity with the property that $Y^n \downarrow Y$:
\begin{align*} Y^n_t = & \phi(t,X_t) + \int_{t}^{T}f(s,X_s,Y^n_s,Z^n_s) + \varphi_n(Y^n_s-\phi(s,X_s))[f(s,X_s,\phi(s,X_s),V_s)+U_s]_-\dd s \\
& - \int_{t}^{T}{Z^n_s \dd W_s}. \end{align*}
Note that it gives an upper approximation for $Y$ while the two previous methods give lower approximations. We can thus hope to define an approximation procedure which gives an upper estimation for the $Y$ (or equivalently for the prices of American options). To keep this property, the time-discretization $(Y^{N},Z^{N})$ is performed in the following way:
\begin{align*}
\begin{cases} \yn_{t_N}= &\phi(t_N,\xn_{t_N}) ,  \\
\yn_{t_k}= & \e_{t_k} \bigg( \yn_{t_{k+1}}+hf(t_k,\xn_{t_k},\yn_{t_{k+1}},\zn_{t_k}) + h\varphi_n(\yn_{t_{k+1}}-\phi(t_{k+1},\xn_{t_{k+1}})) \\
& [ f(t_k,\xn_{t_k},\phi(t_{k+1},\xn_{t_{k+1}}),V^N_{t_k}) + \frac{\phi(t_{k+1},\xn_{t_{k+1}})-\phi(t_k,\xn_{t_k})}{h}]_- \bigg), \\
h \zn_{t_k}= & \e_{t_k}(\yn_{t_{k+1}} \dw^*_k) , \\
 h V^N_{t_k} = & \e_{t_k} \big( \phi(t_{k+1},\xn_{t_{k+1}}) \dw_k^* \big) . \end{cases} \end{align*}
Then, the scheme of Paragraph \ref{par:algo} becomes
\begin{enumerate}
\item  $y^{n,N,R,M}_N(\cdot)=\phi^R(t_N,\cdot)$
\item At a given discretization time $t_k$ ($0 \leq k \leq N-1$), the coefficients $(\alpha_{l,k}^M)$ are defined as the minimizers over $\alpha_l$ of $\displaystyle \frac{1}{M}\sum_{m=1}^{M}{|y_{k+1}^{n,N,R,M}(\xnm_{t_{k+1}})\frac{[\dw_{l,k}^m]_w}{h}-\alpha_l.p_{l,k}^m|^2}.$
Then put $z_{l,k}^{n,N,R,M}(\cdot)=[\alpha_{l,k}^M.p_{l,k}]_z(\cdot)$. 
\item Define $\displaystyle \beta_{l,k}^M= {\rm arginf}_{\beta_l} \frac{1}{M}\sum_{m=1}^{M}{|\phi^R(t_{k+1},\xnm_{t_{k+1}})\frac{[\dw_{l,k}^m]_w}{h}-\beta_l.p_{l,k}^m|^2} .$
\item Finally $\alpha_{0,k}^M$ minimize over $\alpha_0$
\begin{align*}
&\frac{1}{M}\sum_{m=1}^{M}\big| y^{n,N,R,M}_{k+1}(\xnm_{t_{k+1}})+hf^R(t_k,\xnm_{t_k},y_{k+1}^{n,N,R,M}(\xnm_{t_{k+1}}),z^{n,N,R,M}_{l,k}(\xnm_{t_k})) \\
& + h\varphi_n(y^{n,N,R,M}_{k+1}(\xnm_{t_{k+1}})-\phi^R(t_{k+1},\xnm_{t_{k+1}})) [ f^R(t_k,\xnm_{t_k},\phi^R(t_{k+1},\xnm_{t_{k+1}}),[\beta^M_{l,k}.p_{l,k}^m]_v) \\
& + \frac{\phi^R(t_{k+1},\xnm_{t_{k+1}})-\phi^R(t_k,\xnm_{t_k})}{h}]_-  -\alpha_0.p_{0,k}^m\big|^2. \end{align*}
Then put $y^{n,N,R,M}_k(\cdot)=[\alpha_{0,k}^M.p_{0,k}]_y(\cdot)$.
\end{enumerate}
In practice, $\varphi_n$ is a piecewise linear function with $\varphi_n(x)=1$ if $|x| \leq \frac{1}{n}$, $\varphi_n(x)=0$ if $|x| \geq \frac{2}{n}$.

\subsection{Numerical results}
We present here numerical results for American options written on one, three and ten-dimensional assets.
\subsubsection{Dimension $d=1$}
\begin{figure}[h]
\begin{center}
\includegraphics[angle=90,width=0.49\linewidth]{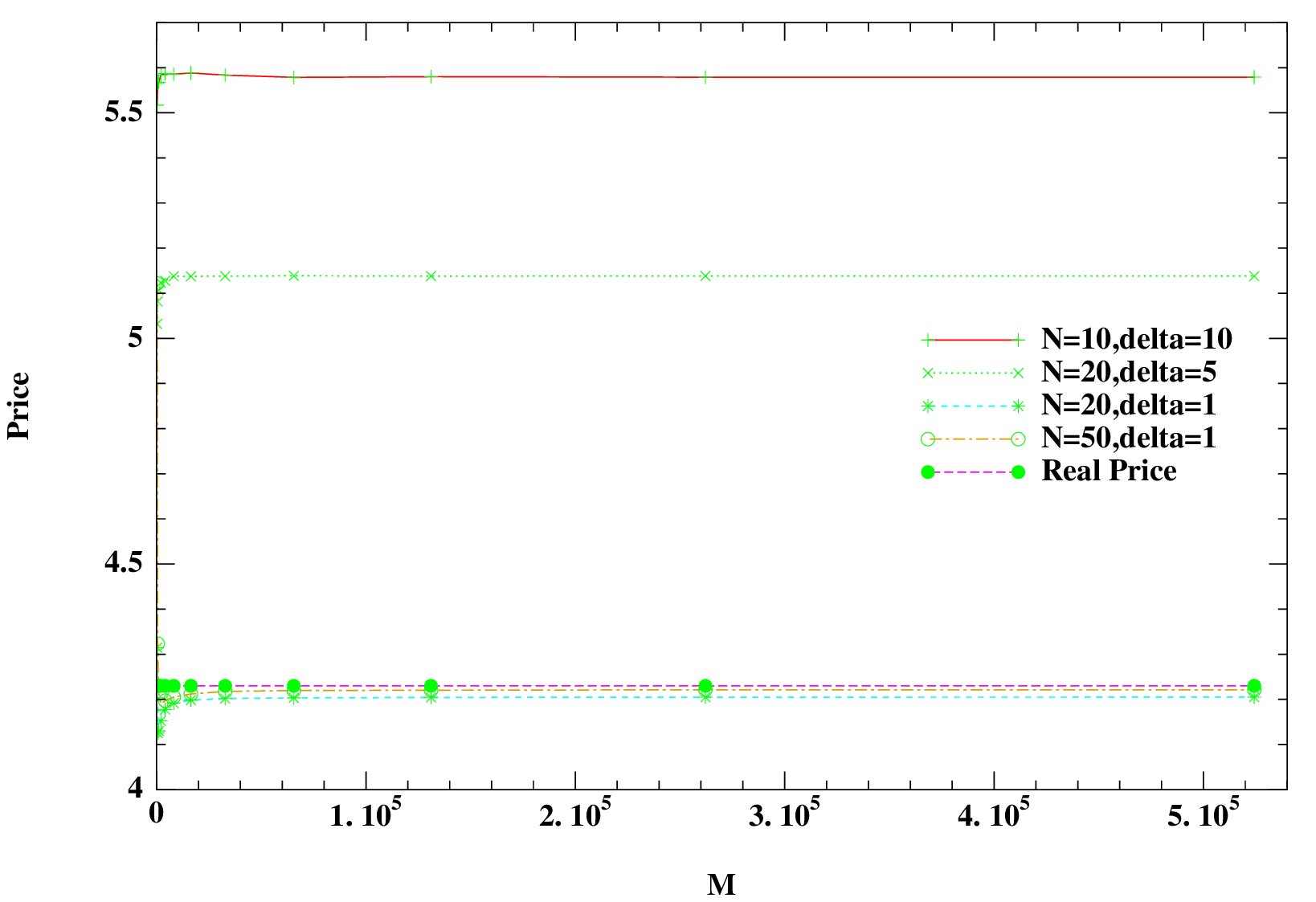}
\includegraphics[angle=90,width=0.49\linewidth]{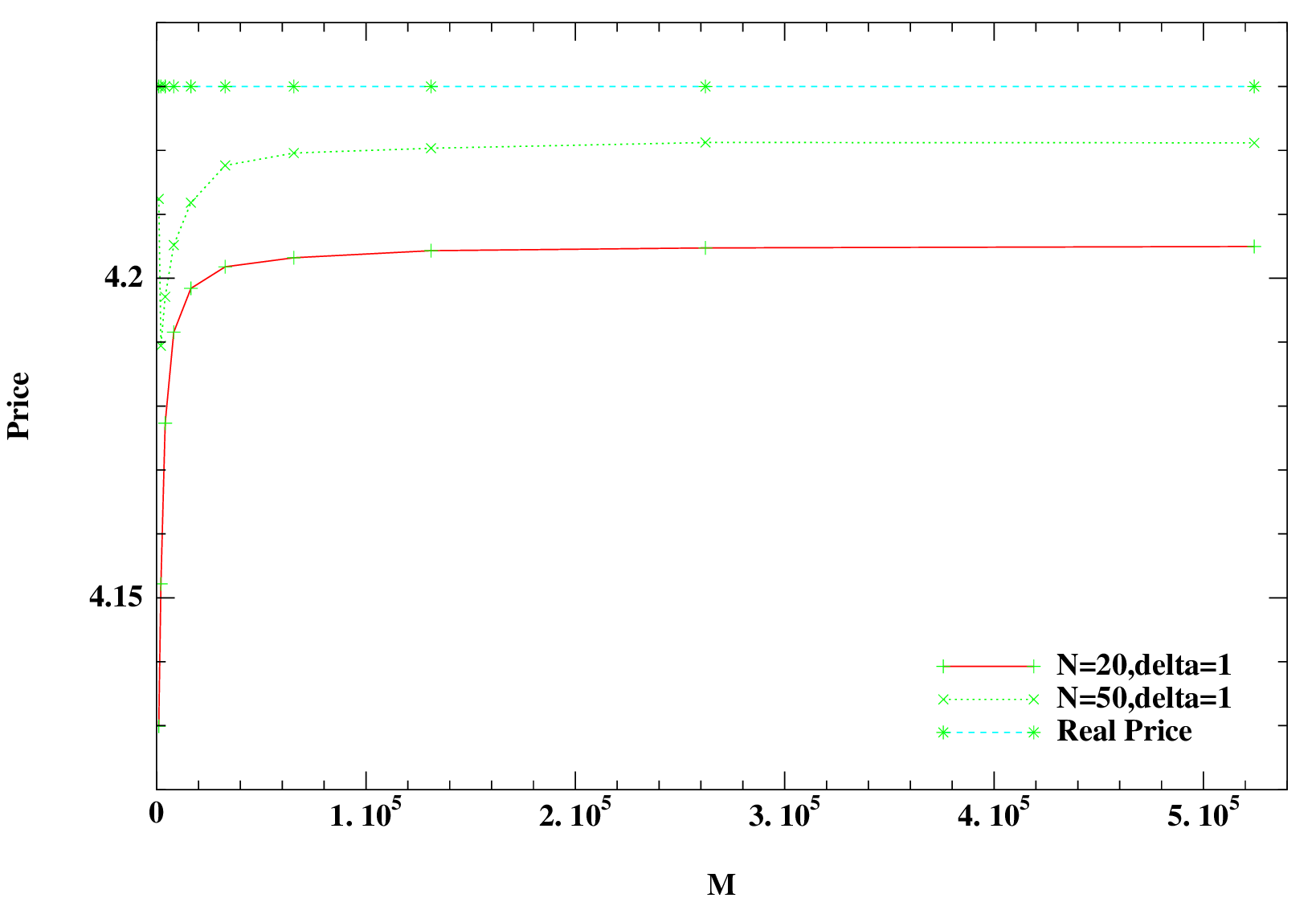}
\end{center}
\caption{American put with $d=1$ \label{fig1}}
\end{figure}

First, we take the same examples than in \cite{bouc:touz:04} where the authors consider the case of a usual Black-Scholes model in dimension $d$ and the payoff $\phi(x)=(K-(\prod_{i=1}^{d}{x_i})^{\frac{1}{d}})_+$. \\
We first take $d=1$, the maturity $T$ of the option is one year, the interest rate $r=0.05$, the volatility $\sigma=0.15$, the strike $K=100$ and $S_0=100$. The reference price obtained by a PDE method is {\bf 4.23} according to \cite{bouc:touz:04} while the price obtained in \cite{bouc:touz:04} by a Malliavin calculus approach to solve the RBSDE is {\bf 4.21}. For our algorithm, we use the basis {\bf HC} of hypercubes with size $\delta$ (see Paragraph \ref{par:HC}). We test the {\em max method}, for different values of $N$ and $\delta$ and we let $M$ increase. On Figure \ref{fig1}, we observe that the bias decreases when $N$ and $\delta$ increase, as predicted in  Theorem \ref{borne erreur}.

\subsubsection{Dimension $d=3$}
\begin{figure}[h]
\begin{center}
\includegraphics[angle=90,width=0.50\linewidth]{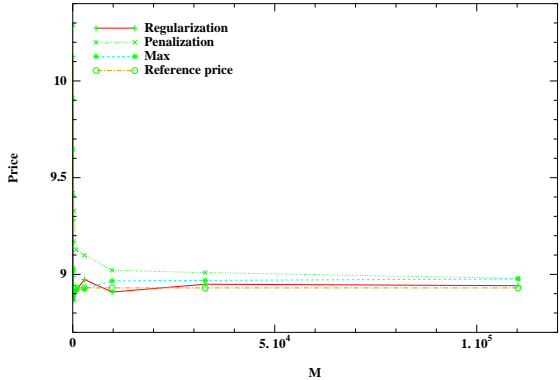}
\end{center}
\caption{American Put with $d=3$ \label{fig2}}
\end{figure}

We deal with the same example with $d=3$. We test the {\em regularization method} with $N=32$, $\delta=9$ and $n=2$, the {\em max method} with $N=44$, $\delta=7$ and the {\em penalization method} with $N=60$, $\delta=2$ and $n=2$. For all of them, we use the basis {\bf HC}. On Figure \ref{fig2}, we observe that the three algorithms give good approximations of the American prices. Nevertheless, we mention that when parameters $(N,\delta,M)$ simultaneously change with the regularization or penalization parameter $n$, the methods may behave very differently and we refer to \cite{lemo:05} for a more complete numerical study.

\subsubsection{Dimension $d=10$}
Finally, we give an example with $d=10=2p$ taken from \cite{ball:page:prin:05}. The model is still a Black-Scholes one : $\frac{dS_t^l}{S_t^l}=(r-\mu_l) \dd t + \sigma_l \dd W_t^l$ for $1 \leq l \leq 10$. The pay-off (or obstacle) is $\max(x_1 \cdots x_p - x_{p+1} \cdots x_{2p},0)$. The interest rate is $r=0$, the dividend rates $\mu_1=-0.05$ and $\mu_l=0$ for $l \geq 2$, the volatility $\sigma_l=\frac{0.2}{\sqrt{d}}$, the maturity is $T=0.5$ and the initial values are $S_0^i=40^{\frac{2}{d}}$, $1 \leq i \leq p$ and $S_0^i=36^{\frac{2}{d}}$ for $p+1 \leq i \leq 2p$. The reference price used in \cite{ball:page:prin:05} is the one derived by \cite{vill:zane:02}, using a PDE method, and is equal to {\bf 4.896}.\\
We use the {\em max method} with a slightly different function basis : instead of just using indicator functions to approximate $Y$, we define also on each hypercube polynomials of degree $1$. With this new function basis, we obtain a price of {\bf 4.876} for $N=60$, $\delta=0.6$ and $M=2^{16}=65536$, within a computational time of $15$ seconds. We note that this algorithm can be quite efficient in high dimension, even with the simplest choice of functions basis.

\section*{Conclusion}
We have proposed a numerical scheme using empirical regressions on function bases. Regarding the model, it is remarkably flexible and allows us to solve generalized BSDEs. Here the stochastic integral in the BSDE is driven by a Brownian motion, but we are preparing an extension of our scheme to deal with general martingales including jumps. Note that the terminal condition of the BSDE is here of the form $\phi(X_T)$ but actually for certain path-dependences, our algorithm may still work after a slight modification (see \cite{gobe:lemo:wari:05}). Our approach is also suitable for reflected BSDEs. However, the robust choice of the {\em penalization} and {\em regularization} parameter as a function of the time step $h$ is a delicate issue, which should deserve a special attention in future works.


\begin{thebibliography}{0}


\bibitem{ball:page:03}
V.~Bally and G.~Pag{\`e}s.
\newblock Error analysis of the optimal quantization algorithm for obstacle problems.
\newblock {\em Stochastic Process. Appl.}, 106(1):1--40, 2003.

\bibitem{ball:02}
V.Bally, M.E. Caballero, B. Fernandez, and N. El Karoui (2002).
\newblock Reflected {BSDE}'s, {PDE}'s and variational inequalities.
\newblock INRIA Report 4455, 2002.

\bibitem{ball:page:prin:05}
V.~Bally, G.~Pag{\`e}s and J.~Printemps (2005).
\newblock A quantization tree method for pricing and hedging multidimensional {A}merican options.
\newblock {\em Mathematical Finance}, 15(1),2005.




\bibitem{barl:buck:pard:97}
G.~Barles, R.~Buckdahn, and E.~Pardoux.
\newblock {Backward stochastic differential equations and integral-partial differential equations}.
\newblock {\em Stochastics and Stochastics Reports}, 60:57--83, 1997.



\bibitem{bouc:touz:04}
B.~Bouchard and N.~Touzi.
\newblock {Discrete time approximation and Monte-Carlo simulation of backward
  stochastic differential equations}.
\newblock {\em Stochastic Processes and their Applications}, 111:175--206, 2004.








\bibitem{elk:97}
N. El Karoui, C. Kapoudjian, E. Pardoux, S. Peng, and M.C. Quenez.
\newblock Reflected solutions of backward {SDE}'s and related obstacle problems for {PDE}'s.
\newblock {\em Ann. Probab. }, 25(2):702--737, 1997.

\bibitem{elka:huan:97}
N.~El~Karoui and S.J.~Huang.
\newblock A general result of existence and uniqueness of backward stochastic differential equations.
\newblock In {\em Backward stochastic differential equations (Paris, 1995--1996)}, Pitman Res. Notes Math. Ser. 364, 27--36, 1997.

\bibitem{elka:peng:quen:97}
N.~El~Karoui, S.G.~Peng and M.C.~Quenez.
\newblock Backward stochastic differential equations in finance.
\newblock {\em Math. Finance}, 7(1):1--71, 1997.


\bibitem{gobe:lemo:wari:05}
E.~Gobet, J.P.~Lemor and X.~Warin.
\newblock A regression-based {M}onte {C}arlo method for backward stochastic differential equations.
\newblock {\em Annals of Applied Probability}, 15(3):2172--2002, 2005.

\bibitem{golu:vanl:96}
G.~Golub and C.F.~Van Loan.
\newblock {\em Matrix computations.3rd ed.}
\newblock {Baltimore, MD: The Johns Hopkins Univ. Press. xxvii, 694 p. }, 1996.

\bibitem{gyor:kohl:krzy:walk:02}
L.~Gy\"orfi, M.~Kohler, A.~Krzyzak and H.~Walk.
\newblock A distribution-free theory of nonparametric regression.
\newblock Springer Series in Statistics, 2002.

\bibitem{hama:jean:05}
S. Hamad\`ene and M. Jeanblanc (2005).
\newblock On the stopping and starting problem : application to reversible investment.
\newblock To appear in {\em Mathematics of Operation Research}, 2005.


\bibitem{jaco:04}
J.~Jacod.
\newblock The Euler scheme for L\'evy driven stochastic differential equations: limit theorems.
\newblock {\em The Annals of Probability}, 32(3):1830--1872, 2004.

\bibitem{lemo:05}
J.P.~Lemor.
\newblock Approximation par projections et simulations Monte-Carlo des \'equations diff\'erentielles stochastiques r\'etrogrades.
\newblock PhD Thesis, Ecole Polytechnique, 2005. http://www.polymedia.polytechnique.fr/Center.cfm?table=These. 

\bibitem{lemo:gobe:wari:05}
J.P.~Lemor, E.~Gobet and X.~Warin.
\newblock Rate of convergence of an empirical regression method for solving generalized backward stochastic differential equations.
\newblock In revision for {\em Bernoulli}.





\bibitem{ma:zhan:05}
\newblock J. Ma and J. Zhang (2005).
\newblock Representations and regularities for solutions to {BSDE}s with reflections.
\newblock {\em Stochastic Processes and Their Applications}, 115(4), 539--569, 2005.



\bibitem{vill:zane:02}
S. Villeneuve and A. Zanette (2002).
Parabolic ADI methods for pricing American option on two stocks.
\newblock {\em Mathematics of Operations Research}, 121--151, Feb 2002. 

\bibitem{zhan:04}
J.~Zhang.
\newblock A numerical scheme for {BSDEs}.
\newblock {\em Ann. Appl. Probab.}, 14(1):459--488, 2004.

\end{thebibliography}
\end{document}